\documentclass[12pt]{article}
\usepackage[latin1]{inputenc}
\usepackage{amsmath,amsthm,amssymb}
\usepackage{amsfonts}
\usepackage{amsmath,amsthm,amssymb,amscd}
\usepackage{latexsym}
\usepackage{color}
\usepackage{grffile}
\usepackage[shortlabels]{enumitem}
\usepackage{graphicx}
\usepackage{enumitem}
\usepackage{mathrsfs}
\usepackage{mathtools}
\usepackage{blindtext}

\makeatletter \@addtoreset{equation}{section} \makeatother

\newtheorem{theorem}{Theorem}[section]

\begin{document}	
	
	
	
		
		
		
		\title{\bf A quick sneak-peek at the $s$-fractional Laplacian operator}
		\author{\bf Debajyoti Choudhuri$^{a}$,
	\\	\small{$^{a}$Department of Mathematics, National Institute of Technology Rourkela, Rourkela - 769008, Odisha, India.}\\
	\small{\it Email: dc.iit12@gmail.com $^a$}}
	 \date{\today}
		 \maketitle
			\begin{abstract}
			\noindent The short note here is to give a few heuristic arguments on the weird looking fractional Laplacian operator. This is certainly going to expand the vision of a reader who is looking to develope a taste for research in this direction.
			\begin{flushleft}
			\end{flushleft}
		\end{abstract}
	\section{Introduction}
	Fractional Laplacian is a fast emerging area of research in many areas of science, mathematics being a major contributor to it owing to the rich theoretical aspects it gives rise to. We will discuss about an application of these operators in {\it atom dislocation in crystals}. One essential reading for this operator is due to {\sc Valdinoci et al} \cite{valdi1}. This short article has been written with minimal details so that the surprising aspects of the $s$-Laplacian operator can be appreciated by the reader(s) without delving into the rigours of the proofs. With my own experience, it is always important for the readers to have a heuristic argument for a better insight for one's understanding.
	\section{Definition}
Following is the definition of the fractional Laplacian:
\begin{align}\label{frac_lap}
\begin{split}
(-\Delta)^su(x)&:=C_{N,s}\int_{\mathbb{R}^N}\frac{u(x)-u(y)}{|x-y|^{N+2s}}dy.
\end{split}
\end{align}
Here, $0<s<1$, $N\geq 1$. The integral, although is singular in nature, is computed in the principle value sense as follows:
\begin{align}\label{frac_lap_PV_sense}
\begin{split}
(-\Delta)^su(x)&:=\underset{\epsilon\searrow 0}{\lim}\int_{\mathbb{R}^N\setminus B_{\epsilon}(x)}\frac{u(x)-u(y)}{|x-y|^{N+2s}}dy.
\end{split}
\end{align}
The constant $C_{N,s}$ in \eqref{frac_lap} is equal to $$C_{N,s}=\frac{4^s\Gamma(\frac{N}{2}+s)}{\pi^{N/2}\Gamma(1-s)}.$$
\subsection{A nice application of $(-\Delta)^s$}
This subsection will introduce the readers to an application of the operator in \eqref{frac_lap}. Consider a lattice of atoms as in the Figure $1$.
\begin{figure}[h]
	\caption{The curve below is the `well of energy'}
	\centering
	\includegraphics[width=0.5\textwidth]{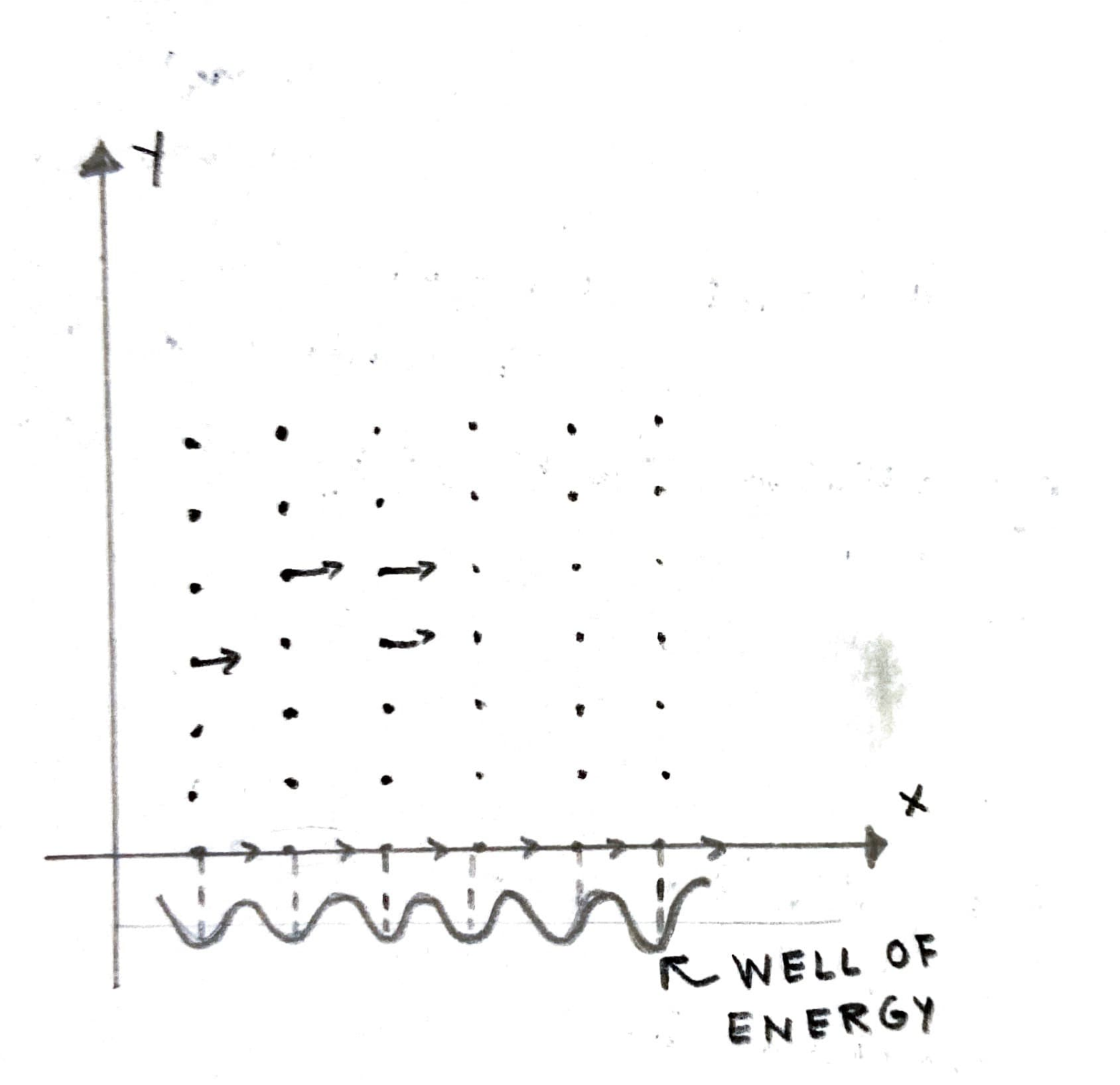}
\end{figure}
We ignore the vertical component of the movement of the atoms. The atoms are distributed in $\mathbb{Z}\times\mathbb{Z}_+$. We denote $V(x,y)$ to be the displacement of the atom in the horrizontal direction. The total energy of the system is then given by
\begin{align}\label{tot_energy}
\mathcal{E}(V)&=\frac{1}{2}\int_{\mathbb{R}_+^2}|\nabla V|^2dxdy+\int_{\mathbb{R}}W(v(y))dy.
\end{align} 
Here $\mathbb{R}_+^2=\mathbb{R}\times(0,\infty)$ and $V(0,y)=v(y)$. Since $\mathcal{E}$ is a $C^1$ functional, therefore, the Fr\'{e}chet derivative of $\mathcal{E}$ is as follows:
\begin{align}\label{tot_energy_der}
\langle\mathcal{E}'(V),\phi\rangle&=\int_{\mathbb{R}_+^2}\nabla V\cdot\nabla \Phi dxdy+\int_{\mathbb{R}}W'(v(x))\phi(x)dy~\text{for every}~\Phi\in H^{1,2}(\mathbb{R}_+^2).
\end{align}
Here $\Phi(0,y)=\phi(y)$. On choosing $\Phi(0,y)=0$, we get 
\begin{align}\label{tot_energy_der'}
0&=\int_{\mathbb{R}_+^2}\nabla V\cdot\nabla \Phi dxdy~\text{for every}~\Phi\in H_0^{1,2}(\mathbb{R}_+^2)\subset H^{1,2}(\mathbb{R}_+^2).
\end{align}
Hence $$-\Delta V=0.$$
Let us now find the critical points of $\mathcal{E}$, i.e. $\langle\mathcal{E}'(V),\phi\rangle=0$. Hence for any $\Phi\in H^{1,2}(\mathbb{R}_+^2)$ we have that 
\begin{align}\label{tot_energy_der_prob}
0&=\int_{\mathbb{R}_+^2}div(\Phi\nabla V)dxdy+\int_{\mathbb{R}}W'(v(x))\phi(x)dy~\text{for every}~\Phi\in H_0^{1,2}(\mathbb{R}_+^2)\\
=&\int_{\mathbb{R}}\phi\frac{\partial V}{\partial y}dx+\int_{\mathbb{R}}W'(v(x))\phi(x)dx.
\end{align}
Therefore, the problem boils down to
\begin{align}\label{boiled_down_prob}
\begin{split}
-\Delta V&=0~\text{in}~\mathbb{R}_+^2\\
\frac{\partial V}{\partial y}(x,0)&=-W'(v(x)).
\end{split}
\end{align}
By the Fourier transform, the problem in \eqref{boiled_down_prob} is equivalent to 
\begin{align}\label{equivalent_boiled_down_prob}
\begin{split}
(-\Delta)^{1/2} V&=W'(v)~\text{in}~\mathbb{R}.
\end{split}
\end{align}
\section{Surprising aspects of solution to \eqref{equivalent_boiled_down_prob}}
We know from the solutions of harmonic functions with Neumann or Dirichlet condition has an exponential decay near infinity. However, we will see that the $s$-harmonic functions have a polynomial decay.\\
{\it Surprise $1$}:~We state the following theorem.
\begin{theorem}\label{asyp_res1}
Suppose $u$ be a solution of 
\begin{align}\label{allen-cahn-type}(-\Delta)^{s} v=W'(v)\end{align}
in $\mathbb{R}$ such that $\underset{x\to-\infty}{\lim}u(x)=0$ and $\underset{x\to\infty}{\lim}u(x)=1$ and $u(0)=1/2$. Then $$\left|u(x)-H(x)+\frac{x}{2sW''(0)|x|^{1+2s}}\right|\leq\frac{C}{|x|^{\theta+2s}}$$
as $|x|\to\infty$. Here $C, \theta>0$ and $H$ is the Heaviside function.
\end{theorem}
\noindent As remarked in the introduction, we will not play with the rigorous proofs of the theorems and the associated results. However, at every stage we will help the reader understand the results via an example or a heuristic but a strong argument to get the flavour of the results.\\
In the case of the Theorem \ref{asyp_res1}, we demonstrate this by taking $u$ to be a compactly supported smooth, analytic function in $\mathbb{R}$. Then $(-\Delta)^su$ is of the order $1/|x|^{2s}$. To see this, we use the following definition of the $s$-fractional operator.
\begin{align}\label{frac_lap_2}
\begin{split}
(-\Delta)^su(x)&:=C_{N,s}\int_{\mathbb{R}}\frac{u(x+y)+u(x-y)-2u(x)}{|y|^{1+2s}}dy.
\end{split}
\end{align}
On expanding the Taylor series of $u$ around $x$ which is near the origin, we get that $(-\Delta)^su(x)\sim|x|^{2-2s}$, i.e. for $s<1$. However, near infinity using the fact that the numerator being bounded, we conclude that the function $1/|y|^{1+2s}$ is integrable if $s>0$. Hence, the operator $(-\Delta)^su(x)$ behaves like $1/|x|^{1+2s}$ near infinity for $0<s<1$.\\
Define $u':=w$. Hence near infinity, the given equation in \eqref{allen-cahn-type} we get 
\begin{align}\label{allen-cahn-type1}(-\Delta)^{s} w=W''(1)w\end{align}
or 
\begin{align}\label{allen-cahn-type2}(-\Delta)^{s} w=W''(0)w\end{align}
In either case we conclude that $w\sim 1/|x|^{1+2s}$. Hence $u\sim 1/|x|^{2s}$. However, this is a contradiction since a compactly supported smooth and analytic function cannot have a polynomial like decay. Therefore, the gist of the Theorem \ref{asyp_res1} is that a solution to \eqref{allen-cahn-type} can neither be compactly supported nor can be a Schwartz class function. It has to decay like a polynomial. This is again quite different from a harmonic function, which decays like an exponential function.\\
{\it Surprise $2$}:~It is well known that the Heat equation decays in the order of exponential function at infinity. Expecting the same to happen for a fractional heat equation will surely take by surprise. Here is a quick look at it. Consider the equation 
\begin{align}\label{heat}
\begin{split}
\frac{\partial}{\partial t}u&=-(-\Delta)^{1/2}u~\text{in}~\mathbb{R}\times(0,\infty)\\
u(x,0)&=\delta_0(x).
\end{split}
\end{align}
As we do it for the heat equation in the local case, we employ the Fourier transform and find the quantity, say at $t=1$, to be $u(x,1)=\mathfrak{F}^{-1}(e^{-|\xi|})=\int_{\mathbb{R}}e^{i\xi\cdot x}e^{|\xi|}d\xi=\frac{2}{1+|x|^2}$. Thus, the decay is like a polynomial again and not like an exponential function.\\
{\it Surprise $3$}:~The problem 
\begin{align}\label{dirchlet_prob}
\begin{split}
(-\Delta)^{s}u&=0~\text{in}~\Omega\\
u&=g~\text{in}~\partial\Omega
\end{split}
\end{align}
is ill-posed. This is not the case with a typical Dirichlet problem that involves a Laplacian. The ill-posedness of the problem is evident from the definition as
\begin{align}
(-\Delta)^su(x)&=\int_{\Omega}\frac{u(x)-u(y)}{|x-y|^{N+2s}}dy+\int_{\mathbb{R}^N\setminus\Omega}\frac{u(x)-u(y)}{|x-y|^{N+2s}}dy~\text{for}~x\in\Omega.
\end{align} 
It is easy to see that $u$ can be chosen to be eccentrically different functions in $\mathbb{R}^N\setminus\Omega$ by keeping $g$ to be the same on the boundary $\partial\Omega$. These different possible extensions of $u$ may lead to different values of $(-\Delta)^su(x)$ for $x\in\Omega$, for an unchanged $g$. This is a clear sign of ill-posedness as $(-\Delta)^su$ becomes multi-valued and doesn't even qualify to be called a function!. Therefore, the correct formulation is
\begin{align}\label{dirchlet_prob1}
\begin{split}
(-\Delta)^{s}u&=0~\text{in}~\Omega\\
u&=g~\text{in}~\mathbb{R}^N\setminus\Omega
\end{split}
\end{align}
or 
\begin{align}
\int_{\Omega}\frac{u(x)-u(y)}{|x-y|^{N+2s}}dy&=:(-\Delta)_{\Omega}^{s}u=0~\text{in}~\Omega\\
u&=g~\text{in}~\partial\Omega.
\end{align} 
{\it Surprise $4$}:~The {\it Harnack's inequality} fails to hold. Consider the following problem.
\begin{align}\label{dirchlet_prob2}
\begin{split}
(-\Delta)^{1/2}u&=0~\text{in}~\mathbb{R}.
\end{split}
\end{align}
Incidentally, 
\[u(x)= \begin{cases}\label{soln_1}
(1-|x|^2)^{-1/2}, & \text{if}~|x|<1 \\
0, & \text{otherwise}
\end{cases}\]
is a solution to \eqref{dirchlet_prob2}. This is because of the following reason. Consider the complex valued function $w(z)=\text{Re}(1-z^2)^{-1/2}$. If one observes carefully the equations \eqref{boiled_down_prob}-\eqref{equivalent_boiled_down_prob}, then it can bee seen that a problem with a $s$-Laplacian can be studied by extending by a dimension and convertig it to a problem with a local operator. For more details on this the reader may refer to {\sc Caffarelli-Silvestre} \cite{caff1}. Therefore,
\begin{align}\label{soln_1_proof}
(-\Delta)^{1/2}u(x)&=-\underset{y\to 0}{\lim}\frac{\partial}{\partial y}w(x+iy)=0.
\end{align}
Note that this solution can never obey the Harnack's inequality.\\
{\it Surprise 5}:~The boundary growth of a $s$-harmonic function is not linear which is unlike a harmonic function. Consider the problem
\begin{align}\label{prob4}
\begin{split}
(-\Delta)^{1/2}u&=1~\text{in}~(-1,1)\\
u&=0~\text{in}~\mathbb{R}\setminus(-1,1).
\end{split}
\end{align}
A solution to \eqref{prob4} is 
\[u(x)= \begin{cases}\label{soln_2}
(1-|x|^2)^{1/2}, & \text{if}~|x|<1 \\
0, & \text{otherwise}.
\end{cases}\]
As like in {\it Surprise $4$}, we choose $w(z)=\text{Re}[(1-z^2)^{1/2}-iz]$. Observe that when $z=iy$, then $w(iy)=(1+y^2)^{1/2}\sim y$ near infinity. Therefore we have to deduct $y$ so that at infinity the solution decays. This is the reason why $iz$ has been deducted while defining the function $w$ (which is not the case with the solution of \eqref{dirchlet_prob2}). Thus
\begin{align}\label{soln_2_proof}
(-\Delta)^{1/2}u(x)&=-\underset{y\to 0}{\lim}\frac{\partial}{\partial y}w(x+iy)=1.
\end{align}
Note that the solution $u$ does not decay linearly near the boundary.
\section{Acknowledgement}
The author acknowledge the facilities received from the Department of Mathematics, NIT Rourkela, India. The author also thanks Prof. Valdinoci for his excellent lectures.
 

\begin{thebibliography}{10}
\bibitem{caff1} L. Caffarelli, L. Silvestre, An extension problem related to the fractional Laplacian, Communications in Partial Differential Equations, 32,  1245--1260, 2007.

\bibitem{valdi1} E. Di Nezza, G. Palatucci, E. Valdinoci, Hitchhiker's guide to the fractional Sobolev spaces, Bull. Sci. Math., 136(5), 521--573, 2012.
 
 
	\end{thebibliography}
\end{document}